\documentclass{amsart}

\usepackage{amsthm}

\usepackage[colorlinks=true,urlcolor=blue,linkcolor=blue,citecolor=blue,breaklinks,pdfusetitle]{hyperref}
\usepackage[hyphenbreaks,anythingbreaks]{breakurl}

\theoremstyle{definition}

\theoremstyle{remark}

\numberwithin{equation}{section}

\begin{document}

\title{Mathematics and the formal turn}

\author{Jeremy Avigad}
\address{Department of Philosophy and Department of Mathematical Sciences, Baker Hall 161, Carnegie Mellon University, Pittsburgh, PA 15213, USA}
\curraddr{}
\email{avigad@cmu.edu}


\date{\today}


\begin{abstract}
Since the early twentieth century, it has been understood that mathematical definitions and proofs can be represented in formal systems systems with precise grammars and rules of use. Building on such foundations, computational proof assistants now make it possible to encode mathematical knowledge in digital form. This article enumerates some of the ways that these and related technologies can help us do mathematics.
\end{abstract}

\maketitle

\section*{Introduction}

One of the most striking contributions of modern logic is its demonstration that mathematical definitions and proofs can be represented in formal axiomatic systems. Among the earliest were Zermelo's axiomatization of set theory, which was introduced in 1908, and the system of ramified type theory, which was presented by Russell and Whitehead in the first volume of \emph{Principia Mathematica} in 1911. These were so successful that Kurt G\"odel began his famous 1931 paper on the incompleteness theorems with the observation that ``in them all methods of proof used today in mathematics are formalized, that is, reduced to a few axioms and rules of inference.'' Cast in this light, G\"odel's results are unnerving: no matter what mathematical methods we subscribe to now or at any point in the future, there will always be mathematical questions, even ones about the integers, that cannot be settled on that basis---unless the methods are in fact inconsistent. But the positive message conveyed by the introduction to the paper is almost as striking: incomplete though they may be, the methods we use can nonetheless be formalized.

The development of the first computational proof assistants around 1970 inaugurated a new era in which it has become possible to carry out formalization in practice. We can now write mathematical definitions, theorems, and proofs in idealized languages, like programming languages, that computers can interpret and check. Young mathematicians everywhere are enthusiastically contributing to the development of digital repositories of mathematical knowledge that can be processed, verified, shared, and searched by mechanical means. And in the mathematical community at large, there is a growing awareness that something important is happening.

This special issue asks: ``Will machines change mathematics?'' From the printing press to to the internet, technology has had major effects on the daily practice of the subject, and it will certainly continue to do so. But it is not hard to justify a negative answer to the question as well: mathematics has always been about coming up with abstractions that enable us think more efficiently, communicate precisely, solve hard problems, and reach a stable consensus as to whether our claims are justified. It's hard to imagine that this will change any time soon. We can use technology to do lots of things, but if we are not using it to do \emph{those} things, then we are probably not doing mathematics.

This essay seeks a middle ground by exploring some of the concrete ways that digital technology can contribute to the mathematical enterprise. I will mention some recent work in the field, but space constraints do not allow a thorough or balanced treatment, so please keep in mind that these are nothing more than scattered examples. The narrative here draws on talks I have given, including the inauguration of the Hoskinson Center at Carnegie Mellon in September 2021 and a presentation to the Topos Institute the following November. But it also draws on talks by many other mathematicians who have recounted some of the reasons that formalization is important to them. See, for example, Thomas Hales's ``Big Conjectures'' at the Big Proof meeting at the Isaac Newton Institute in June 2017; S\'ebastien Gou\"ezel's ``On a Mathematician's Attempts to Formalize his Own Research in Proof Assistants'' at the \emph{Lean Together} meeting at Carnegie Mellon in January 2020; Patrick Massot's ``Why Explain Mathematics to Computers?'' at the Institute Hautes \'Etudes Scientifique in May 2022; Kevin Buzzard's ``The Rise of Formalism in Mathematics'' at the International Congress of Mathematicians in June 2022; Johan Commelin's ``Abstract Formalities'' at the \emph{Fields Medal Symposium in Honor of Akshay Venkatesh}, held at the Fields Institute in November 2022; Adam Topaz's ``The Liquid Tensor Experiment'' at the \emph{Machine Assisted Proof} meeting at the Institute for Pure and Applied Mathematics in February 2023; and Heather Macbeth's ``Algorithm and Abstraction in Formal Mathematics'' at the same meeting. At the time of writing, these are all available on the web and easily found by search engines.

For overviews of interactive theorem proving, see, for example, \cite{avigad:harrison:14,buzzard:20,hales:15}. For a detailed history, see \cite{harrison:urban:wiedijk:14}. The upcoming \emph{Handbook of Proof Assistants}, edited by Jasmin Blanchette and Assia Mahboubi, promises to be a definitive reference.

\section{Verifying theorems}

Proof assistants are used by computer scientists to verify that hardware and software meet their design specifications, a practice known as \emph{formal verification}. The mathematical analogue is verifying the correctness of a mathematical proof, a task which is most compelling when the correctness of the original theorem is in doubt. An early landmark was Thomas Hales's verification of his proof of the Kepler conjecture \cite{hales:05}, which states the the optimal density of a packing of equally-sized spheres in space is attained by the familiar cubic and hexagonal close packing arrangements. Hales launched the project in response to difficulties that arose in the ordinary mathematical review process, some of them stemming from the extensive use of computation in his proof \cite{hales:et:al:17}.

A more recent landmark was the \emph{Liquid Tensor Experiment}, which arose from a challenge posed by Peter Scholze (described in \cite{castelvecchi:21}).
On June 5, 2021, Scholze acknowledged the achievement:\footnote{\burl{https://xenaproject.wordpress.com/2021/06/05/half-a-year-of-the-liquid-tensor-experiment-amazing-developments/}}
\begin{quote}
Exactly half a year ago I wrote the Liquid Tensor Experiment blog post, challenging the formalization of a difficult foundational theorem from my Analytic Geometry lecture notes on joint work with Dustin Clausen. While this challenge has not been completed yet, I am excited to announce that the Experiment has verified the entire part of the argument that I was unsure about. I find it absolutely insane that interactive proof assistants are now at the level that within a very reasonable time span they can formally verify difficult original research.
\end{quote}
For those of us involved in the formalization of mathematics, Scholze's blog post was a major event, a source of encouragement and validation. The challenge was completed on July 14, 2022.

Yet another notable landmark was achieved in 2022, after Thomas Bloom solved a problem posed by Paul Erd\H{o}s and Ronald Graham, showing that any set of natural numbers with positive density has a finite subset whose reciprocals sum exactly to one. Within in a few months, Bloom and Bhavik Mehta verified the correctness of the proof in Lean, while the paper was still under review. In response to a tweet by Kevin Buzzard announcing the result, Timothy Gowers wrote:
\begin{quote}
  Very excited that Thomas Bloom and Bhavik Mehta have done this. I think it's the first time that a serious contemporary result in ``mainstream'' mathematics doesn't have to be checked by a referee, because it has been checked formally. Maybe the sign of things to come?
\end{quote}
As I write this article, Bloom's article has yet to appear in a journal.

Of course, one needs to have confidence that the formal statement that has been verified is an adequate representation of the informal theorem we have in mind. When the statement of a theorem is elementary, as in the case of Bloom's, that is generally not a big concern, but when the statement builds on complex machinery, a more careful audit is called for. For a discussion of the methods used to audit the Liquid Tensor experiment, see the talks by Commelin and Topaz that were mentioned in the introduction, and for discussion of the methods used to audit the formal proof of the Kepler conjecture, see \cite[Section 10]{hales:et:al:17}.

\section{Correcting mistakes}

In addition to helping us verify correctness at the outset, formal verification also helps us discover and correct errors in the mathematical literature. The formalization of the Kepler conjecture led to one major correction \cite{hales:et:al:10}. Similarly, S\'ebastien Gou\"ezel's attempt to formalize a result of Vladimir Shchur on optimal bounds for the Morse lemma in Gromov-Hyperbolic spaces led to the discovery of an error and a nontrivial repair \cite{gouezel:shchur:19}.

The mathematical literature is filled with minor errors and misstatements, and formalization turns up those as well. Hales et al. remark:
\begin{quote}
The detection and correction of
small errors is a routine part of any formalization project. Overall, hundreds of small errors
in the proof of the Kepler conjecture were corrected during formalization. \cite[Section 4.3]{hales:et:al:17}
\end{quote}
Mathematicians have a high tolerance for small errors, but their presence in the literature is a constant source of confusion and frustration. Formalization can help us get the details right.

\section{Gaining insight}

Formalization provides a high degree of confidence in the correctness of a result, but it requires a lot of work, and checking the gritty details of a mathematical argument is not the glamorous side of mathematics. Prospective opportunities for exploration and discovery are much more appealing. A common theme in the literature is that the process of formalization often itself leads to new insights. In his formalization of the four color theorem \cite{gonthier:08}, George Gonthier found it useful to characterize planarity in terms of a combinatorial structure known as a \emph{hypermap}. Hales found the concept useful for reformulating parts of the proof of the Kepler conjecture \cite[Chapter 4]{hales:12}. Methods that emerged from the formalization of the Kepler conjecture, in turn, enabled Hales to prove other longstanding conjectures in discrete geometry, such as the strong dodecahedral theorem \cite[Section 8.6]{hales:12}.

For a more recent example, in his blog post at the midpoint of the Liquid Tensor experiment, Scholze asked himself, rhetorically, ``Did you learn anything about mathematics during the formalization?''
\begin{quote}
  Answer: Yes! The first is a beautiful realization of Johan Commelin. Basically, the computation of the Ext-groups in the Liquid Tensor Experiment is done via a certain non-explicit resolution known as a Breen-Deligne resolution\ldots\,. The Breen-Deligne resolution has certain beautiful structural properties, but is not explicit, and the existence relies on some facts from stable homotopy theory. In order to formalize Theorem 9.4, the Breen-Deligne resolution was axiomatized, formalizing only the key structural properties used for the proof. What Johan realized is that one can actually give a nice and completely explicit object satisfying those axioms, and this is good enough for all the intended applications. This makes the rest of the proof of the Liquid Tensor Experiment considerably more explicit and more elementary, removing any use of stable homotopy theory. I expect that Commelin’s complex may become a standard tool in the coming years.
\end{quote}
He continued to ask himself, ``What else did you learn?''
\begin{quote}
  Answer: What actually makes the proof work! When I wrote the blog post half a year ago, I did not understand why the argument worked, and why we had to move from the reals to a certain ring of arithmetic Laurent series. But during the formalization, a significant amount of convex geometry had to be formalized (in order to prove a well-known lemma known as Gordan’s lemma), and this made me realize that actually the key thing happening is a reduction from a non-convex problem over the reals to a convex problem over the integers.
\end{quote}
Although Scholze did not carry out the formalization himself, it is notable that the project helped him gain a better understanding of his own results. One can argue that this has nothing to do with digital technology per se; spending sufficient time with any piece of mathematics and working with the ideas can help us understand it better. But formal verification generally \emph{requires} us to revisit and revise our proofs. For some,  the appeal of having a pristine, logically complete, and fully-checked artifact at the end is what provides the motivation for putting in the extra work. For others, it is the process itself, and the pleasure of seeing the pieces of a complex mathematical result fit together in perfect harmony.

\section{Building libraries}

Before one can formalize cutting-edge results, one has to formalize all the parts of mathematics they depend on, from elementary arithmetic and algebra to undergraduate and graduate-level mathematics and beyond. Although big-name theorems provide useful milestones, the formal libraries that support them are a better measure of progress. The most impressive mathematical libraries today include those of Mizar, Isabelle, HOL Light, Coq, Metamath, and Lean. By now, these cover most undergraduate-level mathematics and a fair amount of graduate-level mathematics as well. (Efforts to support interoperability between libraries have met with limited success, but it is common practice to manually port results from one system to another.) Lean's library, mathlib, is notable for its open-ended, collaborative model; see \cite{mathlib:20} and the Lean community web pages.

Developing a mathematical library leads to standardization of definitions, notation, and manner of expression, in a manner similar to the way that \LaTeX{} has led to standardization of mathematical notation. We are still learning how to manage the process, and too much uniformity may not be a good thing. But formal repositories provide comprehensive sources of knowledge, helping to turn the cloud into a modern-day library of Alexandria.

\section{Searching for definitions and theorems}

Mathematical results can be expressed in many equivalent forms, and today's search engines are not equipped to detect whether or not one mathematical statement is a straightforward consequence of another. Contemporary \emph{sledgehammer} technology \cite{desharnais:et:al:22}, which combines search techniques with automated reasoning tools, are an important step toward helping us find the definitions and theorems we need. We are also beginning to see evidence that large language models like those behind ChatGPT, when combined with mechanisms for formal oversight, can help; see Section~\ref{section:enabling:ai}.

\section{Refactoring proofs}

Once we have mathematics in a digital format, it can be searched, analyzed, and manipulated by the computer. We can modify a definition and then step through our proofs interactively to make the corresponding changes. We can drop or modify a hypothesis to a theorem and see what breaks. We can extract a particular inference in a long proof, turn it into a separate lemma, and set about generalizing it. We can easily re-order parts of a proof, relying on the assistant to make sure that the hypotheses and intermediate results we need do not get lost in the shuffle. In these ways, a proof assistant is like a word processor for mathematics, allowing us to tinker with and improve our mathematics in much the same way that a word processor allows us to tinker with and improve our writing. Formalizers often speak of ``golfing a proof,'' which is to say, re-engineering it so that it becomes shorter and more efficient. Software engineers use the term ``refactoring'' to describe the analogous process of tidying and improving code that is already in production.

We refactor ordinary mathematical proofs all the time. We reorganize them and break them into smaller, more manageable pieces. We do our best to reduce dependencies between the parts. We extract definitions and lemmas that highlight the key ideas. We look for clever ways of managing information and eliminating clutter. We look for ways to weaken the hypotheses or strengthen the conclusion. We often restructure our proofs to give better indications of how and why the hypotheses are needed. Formalization encourages us to polish and improve our proofs in all these ways, and it provides us with new ways of going about that.

\section{Refactoring libraries}

Even a library's most fundamental concepts can be revised, provided we are willing to go through the library and fix everything that breaks downstream. Despite the effort it requires, communities of formalizers often go back and tinker with definitions and theorems in order to provide better support for the mathematics that is built on top of it.

On the Lean Zulip social media channel, Riccardo Brasca has pointed out that the digital environment changes the way knowledge is updated:\footnote{Leonardo de Moura has highlighted this quotation in some of his presentations on Lean.}
\begin{quote}
  An interesting part of the PR [pull request, i.e.~a request to add a new contribution to the library] is that, with the new theorem, the linter will automatically flag all the theorems that can be generalized (for free!), removing the separability assumption. I think in normal math this is very difficult to achieve. If I generalize a 50 year-old paper that assumes $p \ne 2$ to all primes, there is no way I can manually check and maybe generalize all the papers that use the old one.
\end{quote}
Theorems invariably depend on prior results, and when a theorem in a formal library is strengthened or generalized, a proof assistant enables us to see the effects on downstream results right away and modify them accordingly. In contrast to traditional journals, digital repositories for mathematical knowledge are designed to be revised and updated over time. They are not meant to replace traditional publication models or downplay the importance of archival data, but, rather, to supplement the historical record with useful representations of our current understanding.

\section{Engineering concepts}

Refactoring formal libraries and proofs often yields interesting insights. If you want to formalize the concept of a limit, it helps to use the concept of a filter; the fact that there are so many modes of convergence yields a combinatorial explosion of variants, and the notion of a filter is just the right abstraction to capture them all \cite{holzl:et:al:13}. When one carries out algebraic reasoning formally, it is often more convenient to reason about embeddings between structures rather than reason about substructures, since the former avoids having to transfer results across such embeddings. When we think of a vector as an $n$-tuple, it helps to take $n$ to be any finite type of indices instead of the set of numbers $\{1, \ldots, n\}$, since that avoids having to reindex as often \cite{mathlib:20}. Rather than reasoning with a tower of field extensions $F \subseteq E \subseteq K$, it is more convenient to assume that $K$ is both an $F$- and an $E$-algebra, $E$ is an $F$-algebra, and the actions are suitably compatible \cite{baanen:et:al:22}. When treating the sphere as a manifold, it is convenient to put a chart at every point, because that makes the definition more symmetric and easier to work with.

The observations are mathematically straightforward. For example, the use of filters as a unifying concept goes back to Bourbaki, and contemporary mathematicians have grown used to thinking in structural terms. But when formalizing mathematics, it often takes a lot of insight and experimentation to find just the right away of designing definitions and theorems so that they are reusable and fit together well. Formalizers often have the curious feeling that these types of discoveries are somehow both deep and trivial at the same time.

In her 2023 talk at the Machine Assisted Proof meeting, Heather Macbeth observed that formalization encourages us to identify, name, and reuse abstractions that are implicit in the mathematical literature. With this sort of digital conceptual engineering, the line between software engineering and mathematical understanding is blurred. The innovations often do little more than codify common patterns of reasoning, but this is what mathematical abstraction is all about, namely, identifying patterns and making them explicit. Doing so gives us new ways understanding the mathematics we have, and it provides new resources for future work.

\section{Communicating}
\label{section:communication}

New modes of communication have sprung up around formal technologies. Most proof assistants have mailing lists where users can ask questions or launch discussions, and many now have a presence on social media platforms like Slack, Zulip, and Discord. The Lean community's Zulip channel currently has almost 7,000 subscribers, with around 400-500 active in any given two-week period, and between 500 and 1,000 messages posted every day. People ask questions, pose challenges, share news, and tell jokes. Questions are often answered within seconds. Newcomers are welcomed, and more seasoned users are there to help them as they struggle with the new technology. It is always heartwarming to see those new users, six months later, proudly answering questions and helping the next generation.

Social media and formal technology are only loosely coupled: we certainly don't need proof assistants to communicate mathematics, and platforms like MathOverflow play an important role without any formal component. But digital technology helps to focus communication and social interaction. Mathematical understanding is fleeting; the glorious feeling we have when we finally grasp a complex mathematical argument often vanishes as we struggle to reconstruct the details a week or two later. We may finish a long proof but still harbor fleeting doubts as to whether we have really ironed out all the kinks. In contrast, formal verification brings a satisfying sense of closure. Once we have filled in the last step in a formal proof, there are few things more gratifying than posting it online and sharing it with friends, or, even better, adding it to a communal library. By providing fully precise means of expression, digital technology offers new ways to support mathematical discourse.

\section{Collaborating}

The word ``mathematician'' conjures up the image of a reclusive genius working in isolation, but the average number of authors of a journal publication has been rising steadily over the last several decades \cite{grossman:05,richard:sun:21}. The Liquid Tensor Experiment\footnote{\url{https://leanprover-community.github.io/liquid/}} and a recent formalization of the sphere eversion theorem\footnote{\url{https://leanprover-community.github.io/sphere-eversion/}} \cite{van:doorn:et:al:23} provide models for the way formal digital platforms can support such collaborative efforts. In both cases, the formalization was kept in a public online repository that individual project members could update. Using a system designed by Patrick Massot, participants prepared an informal outline of the proof, also kept online, with hyperlinks to formal statements in the repository. On the same site, collaborators could view a graph displaying dependencies between the definitions and theorems, each color-coded with its status. Definitions were listed in boxes and theorems in circles. A green border indicated that the statement of a result had been formalized, and a green background indicated that the proof has been formalized. A blue border around a result indicated that the statement was \emph{ready} to be formalized, in the sense that all the dependencies were in place, whereas a blue background indicated that the proof was ready to be formalized. Version control systems, commonly used for working on software projects interactively, were used to update both the repository and the blueprint. Collaborators were in constant contact, communicating in real time on social media. The formal proof assistant played a key role, providing a common language and making sure all the pieces fit together.

\section{Managing complexity}

Digital technology can even help us collaborate better with ourselves.
Modularizing mathematical data by breaking it up into smaller pieces with clearly delineated interfaces enables us to deal with one piece at a time, focusing on the information that is relevant at any point in a proof and setting aside the information that isn't \cite{avigad:20}. Proof assistants help us switch context between different parts of a proof by telling us explicitly what information and assumptions are in play, and making sure we use definitions and theorems correctly. The ability to hone in on local features of an argument reduces our cognitive burden; see the contribution by Commelin and Topaz to this issue of the \emph{Bulletin}.

\section{Managing the literature}

At \emph{Machine Assisted Proof}, Jordan Ellenberg moderated a panel discussion with Sophie Morel, Emily Riehl, and Akshay Venkatesh, in which the panelists were prompted to reflect on the prospects for some of the new technologies for mathematics. During the question period, I asked the panelists what they felt were the main challenges facing contemporary mathematics, setting the technology aside. Venkatesh voiced concern about the sheer size of the mathematical literature, and the burden of trying to stay on top of the portions that we need. Riehl agreed, and also noted that it is becoming increasingly difficult for editors and referees to handle the volume of journal submissions. Only Morel offered the problem of ensuring correctness as a primary concern.

Volume of content is not a problem in and of itself. The number of cat videos on YouTube is disconcerting, but nobody feels that they have to watch them all. But as mathematicians we feel a constant pressure to keep on top of the results that we need to know and the methods and techniques that we need to understand in order to do the mathematics we want to do. This requires a kind of informational triage, requiring us to process content at just the level of depth we need to decide where to focus our limited time and energy. The editorial process is an important means by which the mathematical community evaluates and assesses new material, and it is especially important to do so in a manner that is reasoned and fair. We also have to keep up with the literature in order to avoid wasting time and energy by proving results that are already known, possibly in a form that makes it hard for us to recognize them. Letting a student fall into this trap can put an end to an otherwise promising career.

Massot has pointed out that conventional mathematical texts empower authors rather than readers, in the sense that the author decides what information to make available. Formal libraries put the reader in control: since every definition and proof has been reduced to axiomatic primitives, all we need is an interface that allows the reader to display the relevant information on demand. At \emph{Machine Assisted Proof}, Massot demonstrated a prototype tool he has developed with Kyle Miller, which automatically translates a formal proof to one that is written in formulaic but otherwise ordinary English. In a browser, one can click on buttons interspersed through the text, to see all the objects and assumptions in play at that point in the proof. More interestingly, one can click on a plus symbol next to any high-level inference, at which point, the display expands to include a justification of the inference. One can then click on the steps of \emph{that} argument to expand them further, or click on a minus symbol to collapse the text again. This is not an isolated effort; numerous researchers in formal methods have experimented with methods of combining informal exposition with formal content.

In a blog post after the meeting, Terence Tao asked whether AI tools might, in the near future, be able to automatically diagram the structure of a paper on arXiv. It is already commonplace to display dependency diagrams for theorems in a formal library, which users can explore with suitable graphical interfaces. Systems like Isabelle and Coq have long supported tools that enable users to write formal proofs and accompanying Latex and HTML documents in the same file, enabling users to focus on formalizations and informal explanations at the same time. Efforts like these highlight the fact that the goal of formalization is not to replace mathematical exposition but, on the contrary, to enhance it.

\section{Teaching}

When teaching mathematics, quick feedback is almost always helpful to students and almost always hard to come by. When we teach students to write proofs, they usually have to wait a week or more to find out whether their homework is correct, by which point the class has moved on and they are already focused on the next assignment. At best, we can hope that they will give the grader's comments a quick glance before they are buried in a folder.

Interactive proof assistants provide immediate feedback. When a proof is incorrect, or, at least, does not provide enough detail for the computer to assess correctness, students know it right way. When the computer can confirm the correctness of an inference, they also know that right away, and it feels good. Each small frustration is a challenge to overcome, and every small success encourages them to keep going. Students often report that time melts away as they become absorbed in the technology.

Proof assistants are especially good at getting students to recognize the logical form of a mathematical statement and keep track of the objects and assumptions at play at each stage of a proof. As instructors, we are often frustrated when students are confused about the difference between ``every'' and ``some,'' or when they talk about an $x$ or a $\delta$ that hasn't yet entered the scope of an argument. When working with a proof assistant, it's much easier for students to see where they are going wrong and how to fix it.

We have to be careful not to expect too much. As Justin Reich observes in his book, \emph{Failure to Disrupt} \cite{reich:20}, educational technology rarely lives up to its initial hype. It is one thing to teach students how to use a proof assistant and another thing entirely to use a proof assistant to teach mathematics, and to do so in such a way that the system is more of a help than a distraction. We still have a lot to learn, but a growing number of instructors are beginning to experiment with the new technology \cite{kerjean:et:al:22}.

\section{Improving access}

Learning mathematics requires consistent mentorship, and it's hard to become involved in the subject professionally without the right pedigree. The enjoyment and benefits of mathematics are closed off from too many people; even Ramanujan had to travel to England to connect to the mathematical mainstream.

In contrast, adolescents often learn how to write computer programs on their own. Kids can find online tutorials and examples, and then they can experiment and see the effects right away. Online communities also give them opportunities to share their code and insights with others. Formal technology makes doing mathematics more like writing computer programs in these respects. There is hope, then, that the technology can make substantial mathematical engagement accessible to a much broader audience. In Section~\ref{section:communication}, I alluded to the role of online mentorship in helping newcomers getting started with proof assistants. Imagine what we could do if we could find a way to extend that kind of interaction to mathematics as a whole, and foster communities of people coming together to help each other learn and explore. Formalization is not a prerequisite for that, but experience suggests that it can play a role.

We have to be careful: there is a dark side to social media, and we have learned the hard way that crowdsourcing the search for truth doesn't always have the desired effects. We have to avoid thinking that progress will be easy. As Reich points out, even online systems for teaching programming tend to favor students in higher socioeconomic classes, with parents and mentors in a better position to offer encouragement and support. But if not a complete solution, formal methods nonetheless offer novel opportunities for progress.

\section{Using mathematical computation}

We have already grown used to symbolic and numeric computation in science, engineering, industry, and finance, and these are playing an increasing role in pure mathematics \cite{cohn:et:al:17,hales:14,helfgott:15}. Formal methods can make the use of such computation more reliable.

One can always write a computer program to carry out a subtle mathematical calculation, run the program, and then make use of the result, say, in a scientific prediction or a mathematical proof. Often, however, it would be nice to have a stronger guarantee that the program does what we think it does. This brings us to the realm of software verification. Formal methods are commonly used at places like Amazon Web Services, Apple, and Facebook to verify the correctness of code. Using these methods to verify mathematical software is far less common, but it is just as pressing, and more interesting to the mathematical community. For example, proof assistants have been used to verify numeric calculations in number theory \cite{mahboubi:et:al:19}, dynamical systems \cite{immler:18}, and sphere packing \cite{hales:et:al:17}, and symbolic calculations involving Groebner bases \cite{pottier:08}, semidefinite programming \cite{harrison:07}, convex optimization \cite{bentkamp:mir:avigad:23}, and algebra \cite{mahboubi:sibut:pinote:21}.

Verification should become the standard for mathematical computation. Why settle for less? As Hales et al.~point out, the ability to trust the results of computation encourages us to offload parts of our mathematical proofs:
\begin{quote}
In the original [proof of the Kepler conjecture], computer calculations were a last resort after as much was
done by hand as feasible. In the blueprint, the use of computer has been
fully embraced. As a result, many laborious lemmas of the original proof
can be automated or eliminated altogether. \cite[Section 4.3]{hales:et:al:17}
\end{quote}
See also Heather Macbeth's talk at \emph{Machine Assisted Proof}.

Sometimes even just \emph{specifying} the intended behavior of a piece of mathematical software in a proof assistant is helpful, since that links the code to a precise mathematical description. The meaning of an expression in a computer algebra system or a function in a library for numeric computation is often left ambiguous by the system's documentation. Without a precise specification, it isn't even meaningful to ask whether the code is correct; in the words of Wolfgang Pauli, it is ``not even wrong.'' Whether or not we ultimately verify a piece of mathematical software, a formal specification can be used to relate the results to proofs, scientific results, and other pieces of software.

\section{Using automated reasoning}
\label{section:enabling:automated:reasoning}

In the field of artificial intelligence, there is a longstanding distinction between \emph{symbolic methods}, which rely on explicit, humanly-interpretable representations of knowledge, and \emph{subsymbolic methods}, in which information is represented by patterns of data that may be distributed across a neural network. Formal methods fall squarely under the first category, since they use precise, logic-based representations and rules of inference. John Haugeland dubbed such approaches \emph{GOFAI}, short for \emph{Good Old Fashioned AI}, in 1985 \cite{haugeland:85}.

Good old fashioned AI has given rise to a range of automated reasoning tools, including first-order theorem provers, SAT solvers, SMT solvers, and term rewriters. These tools address different problem domains, but what they have in common is that queries are posed in symbolic languages with precise, mathematically specified semantics. Moreover, the search for a solution is designed to guarantee its validity, assuming the search is implemented correctly. The main challenge to system engineers is combinatorial explosion in the search space, which one attempts to manage with heuristics. Such automated reasoning tools are commonly used for hardware and software verification, as well as for various problem-solving tasks like planning and scheduling.

To date, however, notable applications of automated reasoning to mathematics have been few and far between \cite{avigad:18,heule:kullmann:17}. One application is the development of \emph{sledgehammer} tools that are used to fill in steps in formal proofs automatically \cite{desharnais:et:al:22}. Reasoning with a proof assistant requires combing the library for theorems and painstakingly spelling out inferences that are often intuitively clear. A sledgehammer uses heuristics to select a manageable set of potentially relevant facts from the library and then asks an automated reasoning tool to spell out the details automatically. The system then tries to reconstruct the justification within the proof assistant, checking every step against the axiomatic foundation.

But the digitization of mathematics opens the door to broader uses of automated reasoning as well. Once we have our mathematics in digital form, we are in a position to express mathematical problems in such a way that we can explore the use of automation to justify more substantial inferences or even discover new results. In the years to come, proof assistants and better interfaces should make these tools more accessible to the mathematical community. Once these tools are in the hands of capable mathematicians, there is no telling what might happen.

Automated reasoning tools invariably have bugs, and proof assistants can also provide means of verifying the correctness of the results. A few years ago, Joshua Brakensiek, Marijn Heule, John Mackey, and David Narv\'aez used a SAT solver to resolve the remaining case of a problem known as Keller's conjecture, a conjecture in discrete geometry that claims that a certain statement holds of $n$-dimensional Euclidean space, for every $n$ \cite{brakensiek:et:al:20}. It turns out that the conjecture fails in dimension 8 and above; Brakensiek et al.~established that the result holds up to dimension 7. Cayden Codel, a PhD student at Carnegie Mellon, is developing a formal library that can be used to verify the correctness of encodings used with SAT solvers, and Joshua Clune has formally verified the pen-and-paper reduction of the geometric statement to a graph theoretic one, which was instrumental to the solution \cite{clune:23}. An end-to-end verification of the Keller problem, including both the conventional mathematical component and the computer-assisted part, seems to be within reach.

\section{Using machine learning}
\label{section:enabling:ai}

Formalization can also serve as a gateway to the use of subsymbolic methods of AI for mathematics. Those of us working in formal methods can't help but view the stunning successes of machine learning over the last decade with a measure of anxiety, wondering what roles will remain for symbolic methods in the years to come. But applications of machine learning to mathematics have had only mixed success so far, and it seems that a number of important tasks will require a synthesis of the machine learning and symbolic approaches.

For example, machine learning has been used to decide which theorems in a large formal library to use with a sledgehammer, a task known as \emph{premise selection}. (See \cite{goertzel:et:al:22,mikula:et:al:23}, for some recent examples.) A large language model like ChatGPT can even be used to draft an outline of a formal proof, leaving a sledgehammer to fill in the details \cite{jiang:et:al:22}. This should be a more effective route to discovering proofs than asking a large language model to write informal mathematics because a formal proof checker provides a clear signal as to whether the result is successful. Two important goals in the use of AI for mathematics are to achieve \emph{autoformalization}, the ability of a system to translate informal mathematics to formal mathematics automatically, and \emph{self-learning}, the ability of a system to explore mathematics on its own and learn from its successes and failures. In both cases, relying on a formal checker to provide feedback and guidance may be key to making progress.

Mathematicians have also begun using machine learning to explore mathematics in other ways, finding patterns in data and searching for mathematical objects. For example, working with Alex Davies and other colleagues at DeepMind, Andr\'as Juh\'asz and Marc Lackenby used a neural network to detect a relationship between certain types of knot invariants, enabling them to formulate and prove a new theorem \cite{davies:et:al:21}. Similarly, Geordie Williamson was able to explore the relationship between certain graphs and polynomials that arise in representation theory \cite{davies:et:al:21}. Adam Zsolt Wagner was able to use machine learning to find counterexamples to a number of conjectures in graph theory \cite{wagner:21}. To a computer, constructing exotic graphs is like playing a game; the mathematician's task is to design the game and let the computer learn to play it.

The examples in the last paragraph do not rely on formalization, but they do require designing computational representations of mathematical objects in such a way that their properties are available to mechanical exploration. It seems reasonable to hope that formalization of mathematical knowledge will help pave the way to designing such representations.

\section{Supporting a synthesis of machine learning and symbolic AI}
\label{section:supporting:a:synthesis}

Because it provides a platform for the synthesis of machine learning and symbolic reasoning, formal methods in mathematics may give rise to benefits that extend far beyond mathematics itself. Machine learning and symbolic AI have complementary strengths: systems for machine learning can process vast amounts of data but are not good at justifying their results and getting the details right, while symbolic AI, while tending to the low-level inferential details, fails to incorporate intuition and insight.

It is interesting to note that in his seminal paper on machine intelligence \cite{turing:50}, Alan Turing highlighted the two ends of the spectrum as good places to begin:
\begin{quote}
We may hope that machines will eventually compete with men in all purely intellectual fields. But which are the best ones to start with? Even this is a difficult decision. Many people think that a very abstract activity, like the playing of chess, would be best. It can also be maintained that it is best to provide the machine with the best sense organs that money can buy, and then teach it to understand and speak English. This process could follow the normal teaching of a child. Things would be pointed out and named, etc. Again I do not know what the right answer is, but I think both approaches should be tried.
\end{quote}
This passage is eerily prescient, given that the success of symbolic methods in beating the World Chess Champion, Garry Kasparov, was a major landmark in the 1990s and the successes of large language models today have again revolutionized the field.

Arguably the most important challenge to artificial intelligence now is to find approaches that combine the best features of machine learning and symbolic reasoning. Mathematics is an ideal place to start: mathematical reasoning starts with ideas, insights, and intuitions that draw on mathematical experience and understanding, and then it subjects them to precise rules of inference and argumentation. This dual nature was a common theme in a recent report of the US National Academies of Sciences, Engineering, and Medicine titled \emph{AI to Assist Mathematical Reasoning} \cite{national:academies:23}.

The final sentence of Turing's article is a poignant today as it was in his time:
\begin{quote}
We can only see a short distance ahead, but we can see plenty there that needs to be done.
\end{quote}
Given the growing role that artificial intelligence is bound to play in our daily lives, it is as important as ever that we appreciate the value of mathematical thought and keep it central to our decision-making processes.

\section*{Conclusions}

Formal technologies for mathematics are still new, and there is no reason to oversell them. Proof assistants are not easy to use, and many of the challenges they pose have little to do with mathematics per se. Formalization can be a pain in the neck, and most cutting-edge mathematics is, at present, out of reach.

But I hope I have convinced you that the technology holds a lot of promise, and that even now it offers a number of benefits. It can help us verify correctness, build libraries, explore concepts, collaborate, teach, compute, and discover new theorems. These have been part of what it means to do mathematics for as long as the subject has been around, and insofar as digital technology can help, it is here to stay.

An important question facing the mathematical community is how best to make progress. Younger generations are surprisingly comfortable with the new technology, and increasing numbers of graduate students and even undergraduate students are gravitating to it. They bring with them enormous amounts of talent, enthusiasm, energy, and creativity, and they are doing impressive things. Many of them will one day look for jobs in mathematics departments, and we need to think about how to value and reward their contributions.

It is misleading to say that technology will change mathematics, because only we can do that. But the availability of the new technology will undoubtedly alter our conception what it means to do mathematics, and we should be mindful of those changes. The fact that the subject has been able to evolve over the centuries while preserving its essence is a testament to its resilience. After all, coming up with new ways of understanding is what mathematics is all about.

\bibliographystyle{amsplain}
\bibliography{formal_turn_arxiv}
\end{document}